\documentclass{article}
\usepackage{amsmath,amssymb}
\usepackage{enumerate}
 \usepackage{amsfonts}
 \usepackage{amssymb}
 \usepackage{amsthm}
 \usepackage{graphicx}

\makeatletter
\def\fps@figure{htbp}
\makeatother
\usepackage[text={28cc,42cc},centering]{geometry}
\newcommand\toe{\stackrel{e}{\to}}
 \newcommand\tov{\stackrel{v}{\to}}

\author{Nikolay Kolev \thanks{Supported by the Scientific Research Fund of the St. Kl.
Ohridski Sofia University under contract 90-2008}}

\begin{document}

\title{NEW UPPER BOUND FOR THE EDGE FOLKMAN NUMBER $F_e(3,5;13)$}

\maketitle

\begin{abstract}
  For a given graph $G$ let $V(G)$ and $E(G)$ denote the vertex and
  the edge set of $G$ respevtively.  The symbol $G \toe (a_1, \ldots ,
  a_r)$ means that in every $r$-coloring of $E(G)$ there exists a
  monochromatic $a_i$-clique of color $i$ for some $i \in
  \{1,...,r\}$.  The edge Folkman numbers are defined by the equality
$$ F_e(a_1, \ldots , a_r ;q) = \min\{ | V(G)| : G \toe  (a_1, \ldots ,
   a_r; q) \text{ and } cl(G)<q  \}.
$$
In this paper we prove a new upper bound on the edge Folkman number
$F_e(3,5;8),$ namely

$$  F_e(3,5;8) \leq  21$$
This improves the bound

$$  F_e(3,5;8) \leq  24,$$
proved in \cite{KNgod}


    \textbf{2000 Math.\ Subject Classification.} 05C55
\end{abstract}

\section{Introduction}

Only finite non-oriented graphs without multiple edges and loops are
considered. We call a $p$-clique of the graph $G$ a set of $p$
vertices each two of which are adjacent. The largest positive
integer $p$ such that $G$ contains a $p$-clique is denoted by cl(G).
A set of vertices of the graph $G$ none two of which are adjacent is
called an independent set. The largest positive integer $p$ such
that $G$ contains an independent set on $p$ vertices is called the
independence number of the graph $G$ and is denoted by $ \alpha(G).$
In this paper we shall also use the following notations:
\begin{itemize}
  \item $V(G)$ is the vertex set of the graph $G$;
  \item $E(G)$ is the edge set of the graph $G$;
  \item $N(v)$, $v\in V(G)$ is the set of all vertices of $G$ adjacent to $v$;
  \item $G[V]$, $V \subseteq V(G)$ is the subgraph of $G$ induced by $V$;

  \item $K_n$ is the complete graph on $n$ vertices;
  \item $\overline{G}$ is  the complementary graph of $G.$
  \end{itemize}

Let $G_1$ and $G_2$ be two graphs without common vertices. We denote
by $G_1+G_2$ the graph $G$ for which $V(G)=V(G_1)\cup V(G_2)$ and
$E(G)=E(G_1)\cup E(G_2)\cup E'$ where $E'=\{[x,y]:x\in V(G_1), y\in
V(G_2)\}$. It is clear that

\begin{equation}
 cl(G_1 + G_2) = cl(G_1) + cl(G_2).
\end{equation}



\medskip

\textbf{ Definition 1.} Let $a_1, \ldots , a_r$ be positive integers.
The symbol $G \tov (a_1, \ldots , a_r)$ means that in every
$r$-coloring of $V(G)$ there is a monochromatic $a_i$-clique in the
$i$-th color for some $i \in \{1, \ldots, r\}$.
\medskip


\medskip

\textbf{ Definition 2.} Let $a_1, \ldots , a_r$ be positive
integers. We say that an $r$-coloring of $E(G)$ is $(a_1, \ldots ,
a_r)$-free if for each $i=1, \dots , r$ there is no monochromatic
$a_i$-clique in the $i$-th color. The symbol $G \toe (a_1, \ldots ,
a_r)$ means that there is no $(a_1, \ldots , a_r)$-free coloring of
$E(G).$

\medskip

 The smallest positive integer $n$ for which $K_n \toe (a_1, \ldots ,
   a_r)$ is  called a Ramsey number and is denoted by $R (a_1, \ldots ,
   a_r).$  Note that the Ramsey number $R(a_1, a_2)$ can be
   interpreted as the smallest positive integer $n$ such that for
   every $n$-vertex graph $G$ either $cl(G)\geq a_1 $ or $\alpha(G)\geq a_2.
   $ The existence of such numbers was proved by Ramsey in \cite{R}. We
   shall use only the values
   $R(3,3)=6$ and $R(3,4)=9,$ \cite{GG}.

   The edge Folkman numbers are defined by the equality
$$ F_e(a_1, \ldots , a_r ;q) = \min\{ | V(G)| : G \toe  (a_1, \ldots ,
   a_r; q) \text{ and } cl(G)<q  \}.
$$
It is clear that $G\toe (a_1,\ldots ,a_r)$ implies $cl(G)\geq
\max\{a_1,\ldots , a_r\}$. There exists a graph $G$ such that $G
\toe (a_1, \ldots , a_r) $ and $cl(G) = \max\{a_1,\ldots , a_r\}$.
In the case $r=2$ this was proved in \cite{F} and in the general
case in \cite{NR}. Therefore
$$
  \label{1} F_e(a_1 ,\ldots , a_r ;q) \mbox{ exists
    if and only if }
  q> \mbox{max} \{ a_1,\ldots , a_r \}.
$$

It follows from the definition of $R(a_1,\dots,a_r)$ that
\[
F_e(a_1,\dots,a_r;q)=R(a_1,\dots,a_r) \text{ if }q>R(a_1,\dots,a_r).
\]

The smaller the value of $q$ in comparison to $R(a_1,\dots,a_r)$ the
more difficult the problem of computing the number
$F_e(a_1,\dots,a_r;q).$

Among the edge Folkman numbers of the kind
$F_e(a_1,\dots,a_r;R(a_1,\dots,a_r))$ only the following ones are
known :
\begin{align*}
&F_e(3,3;6)=8, \text{ }\cite{G}; \\
&F_e(3,4;9)=14, \text{ } \cite{N349};  \\
&F_e(3,5;14)=16,\text{ } \cite{L};  \\
&F_e(4,4;18)=20, \text{ }\cite{L}; \\
&F_e(3,3,3;17)=19 \text{ } \cite{L}.\\
\end{align*}

Only three  edge Folkman numbers of the kind
$F_e(a_1,\dots,a_r;R(a_1,\dots,a_r)-1)$ are known, namely
$F_e(3,4;8)=16,$ $F_e(3,3;5)=15$ and $F_e(3,3,3;16)=21$. The number
$F_e(3,4;8)=16,$ was computed in the papers \cite{KNdokl},
\cite{KNgod}. The inequality $F_e(3,3;5)\le 15$ was proved in
\cite{N335} and the inequality $F_e(3,3;5)\ge 15$ was obtained by
the means of computer in \cite{PRU}. The inequality
$F_e(3,3,3;16)\ge 21$ was proved in \cite{L} and the opposite
inequality $F_e(3,3,3;16)\le 21$ in \cite{N33316}. At the end of
this exposition we shall mention that we know only one edge Folkman
number of the kind $F_e(a_1,\dots,a_r;R(a_1,\dots,a_r)-2),$ namely
$F_e(3,3,3;15)=23$, \cite{N33315} and only one edge Folkman number
of the kind $F_e(a_1,\dots,a_r;R(a_1,\dots,a_r)-3),$ namely
$F_e(3,3,3;14)=25$, \cite {N33314}. No other edge Folkman numbers
are known.

     This paper is dedicated to the Folkman number $F_e(3,5;13).$

The best known lower  bound on this number is   $F_e(3,5;13) \geq
18,$ which was proved by Lin in \cite{L}. Later Nenov proved in
\cite{N3518} that equality $F_e(3,5;13)= 18$ can be achieved  only
for the graph $K_8+C_5+C_5$. Thus if $K_8+C_5+C_5 \toe (3,5)$   then
$F_e(3,5;13)= 18$ and otherwise $F_e(3,5;13) > 18.$ So far  nobody
was able to check whether $K_8+C_5+C_5 \toe (3,5)$. The best known
upper bound is  $F_e(3,5;13) \leq 24,$ \cite{KNgod}.

We consider the graph $Q$, which was first introduced in \cite{GG}
and whose complementary graph is given in the picture below.

 We shall use the following properties of the graph $Q$:
\begin{gather}
cl(Q)=4, \text{ }\cite{GG};\\
 \alpha(Q)=2, \text{ } \cite{GG};\\
 Q \tov (3,4),\text{ } \cite{N34v}.
\end{gather}

The goal of this paper is to prove the following

\medskip

\textbf{Theorem} Let $G = K_8 + Q$. Then $G \toe (3,5)$.

\medskip

It follows from (1) and (2) that $cl(G)=12$. Since $|V(G)| = 21$ we
obtain from the Theorem the following

\medskip

\textbf{Corollary} $F_e(3,5; 13) \leq 21$.

\section{Proof of the theorem}
 Assume the opposite: that there exists a $(3,5)$-free 2-coloring of
$E(G) $
\begin{equation}
 E(G)=E_1\cup E_2, \qquad E_1\cap E_2=\emptyset.
\end{equation}
We shall call the edges in $E_1$   blue and the edges in $E_2$  red.

We define for an arbitrary vertex $v\in V(G)$ and index $ i=1,2 :$
\begin{align*}
N_i(v)&=\{x\in N(v)\mid [v,x]\in E_i\},
 \\
G_i(v)&=G[N_i(v)]\\
A_i(v) &= N_i(v) \cap V(Q)
\end{align*}

Let $H$ be a subgraph of $G$. We say that $H$ is a monochromatic
subgraph in the blue-red coloring (5) if $E(H)\subseteq E_1$ or
$E(H)\subseteq E_2$. If $E(H)\subseteq E_1$ we say that $H$ is a
blue subgraph and if $E(H)\subseteq E_2$ we say that $H$ is a red
subgraph.

It follows from the assumption that the coloring (5) is (3,5)-free
that
 \begin{equation}
      cl(G_1(v)) \leq 4
 \text{ and }
   cl(G_2(v)) \leq 8
   \text{ for each }  v \in V(G)
 \end{equation}
 Indeed, assume that $cl(G_1(v)) \geq 5.$ Then there must be no blue
 edge connecting any two of the vertices in $ cl(G_1(v)) $ because
 otherwise this blue edge together with the vertex $v$ would give a
 blue triangle. As we assumed $cl(G_1(v)) \geq 5$ then we have a red
 5-clique. Analogously assume $cl(G_2(v)) \geq 9$. Since $R(3,4) =9$,
 then we have either a blue 3-clique or a red 4-clique in $G_2(v).$ If
 we have a blue 3-clique in $G_2(v)$ then we are through. If we have
 a red 4-clique   then this 4-clique together with the vertex $v$ gives a
 red 5-clique. Thus (6) is proved.

 We shall prove that
\begin{equation}
 cl(G[A_1(v)])+cl(G[A_2(v)])\leq 5
 \text{ for each }
  v \in V(K_8)
 \end{equation}
Assume that (7) is not true, i.e. that there exists a vertex $v \in
V(K_8)$ such that

$$cl(G[A_1(v)])+cl(G[A_2(v)])\geq 6.$$

Then as there  are seven more vertices in $  V(K_8)$ with the
exception of $v,$ it follows that
$$cl(G_1(v)) +cl(G_2(v)) \geq 13.             $$
It follows from the pigeonhole principle that either $cl(G_1(v))\geq
5$ or $cl(G_2(v))\geq 9,$ which contradicts (6). Thus (7) is proved.

Now we shall prove that
    \begin{equation}
 cl(G[A_1(v)])=4 \text{ or } cl(G[A_2(v)])=4
 \text{ for each }  v \in V(K_8)
 \end{equation}

By (2) we have

\begin{equation}
  cl(G[A_i(v)]) \leq 4,~ i =1,2.
\end{equation}

Assume that (8) is not true. Then we obtain from (9) that

\begin{equation}
  cl(G[A_1(v)])\leq 3 \text{ and } cl(G[A_2(v)]) \leq 3~\text{for some}~v \in  V(K_8).
 \end{equation}

It follows from (4) that in every 2-coloring of $V(Q),$ in which
there are no $4$-cliques in none of the two colors then there are
$3$-cliques in the both colors. Therefore the inequalities in (10)
are in fact equalities, which contradicts (7). Thus (8) is proved.

        We shall prove that there are at least  7   vertices   $  v \in V(K_8)$   such
that \[cl(G[A_2(v)]) = 4.\]
 Assume the opposite. Then it follows from
(8) that  there are at least 2 vertices  $ v_1, v_2$ in $V(K_8)$
such that $ cl(G[A_1(v_1)]) =cl(G[A_1(v_2)]) = 4.$ Now we conclude
  from (6) that all edges from $v_1,$ $v_2$ to all vertices in
$V(K_8)$ (including the edge $[v_1, v_2]$) are red. Since $R(3,3)=6$
there is a monochromatic $3$-clique in the other 6 vertices  in
$V(K_8)$ excluding  $v_1,$ $v_2.$ If this monochromatic $3$-clique
is blue then we are through. If it is red then this monochromatic
$3$-clique together with the edge $[v_1, v_2]$ forms a red
$5$-clique which is a contradiction. Thus we proved that there are
at least 7 vertices $ v \in V(K_8)$   such that $cl(G[A_2(v)]) = 4.$

    We obtain  from $R(3,3)=6$ that there is a monochromatic $3$-clique among these 7 vertices sufficing
     $cl(G[A_2(v)]) = 4.$ This
    $3$-clique is red (otherwise we are through). Let us denote its
vertices by $a_1,$  $a_2,$  $a_3.$ It follows from (7) that

$$ cl(G[A_1(a_i)]) \leq 1,\text{  }  i=1,2,3.$$
Now we have from (3) that $|A_1(a_i) | \leq 2.$ Then there are at
least 7 vertices in $V(Q)$ from which there are only red edges to
$a_1,$ $ a_2,$ $a_3.$ As $R(3,3)=6$ and $\alpha(Q)=2$ it follows
that there is a $3$-clique among these 7 vertices. If this
$3$-clique is monochromatic blue then we are through. Therefore it
is not monochromatic blue and hence there is a red edge in it. This
red edge together with  $a_1,$ $a_2,$  $a_3$ gives a monochromatic
red 5-clique.

\medskip

The Theorem is proved.

\medskip

\textbf{Acknowledgements.} I am grateful to prof. N. Nenov whose
important comments improved the presentation of the paper.

\begin{flushleft}
Faculty of Mathematics and Informatics\\
St.~Kl.~Ohridski University of Sofia\\
5, J.~Bourchier Blvd.\\
BG-1164 Sofia, Bulgaria\\
e-mail: \texttt{nickyxy@fmi.uni-sofia.bg}
\end{flushleft}

\end{document}